\documentstyle[11pt]{article}

\newcommand{\be}{\begin{equation}}
\newcommand{\ee}{\end{equation}}
\newcommand{\bea}{\begin{eqnarray}}
\newcommand{\ea}{\end{eqnarray}}

\newcommand{\ncm}{\newcommand}

\ncm{\rncm}{\renewcommand}

\rncm{\theequation}{\thesection.\arabic{equation}}
\ncm{\setc}[1]{\setcounter{equation}{#1}}
\rncm{\sec}{\setc{0}\section}

\ncm{\beq}{\begin{equation}}

\ncm{\eeq}{\end{equation}}
\ncm{\no}[1]{(\ref{#1})}

\ncm{\beanon}{\begin{eqnarray*}}

\ncm{\eea}{\end{eqnarray}}

\ncm{\eeanon}{\end{eqnarray*}}

\ncm{\qed}{\hspace*{\fill}\rule{0.24cm}{0.24cm}}

\ncm{\R}{{\cal R}_{H\o H}}
\ncm{\G}{{\cal G}}
\ncm{\J}{{\cal J}}
\rncm{\H}{{\cal H}}
\ncm{\A}{{\cal A}}
\ncm{\B}{{\cal B}}
\ncm{\C}{{\cal C}}
\ncm{\D}{{\cal D}}
\ncm{\F}{{\cal F}}
\ncm{\K}{{\cal K}}
\ncm{\T}{{\cal T}}
\ncm{\N}{{\cal N}}
\ncm{\M}{{\cal M}}
\rncm{\L}{{\cal L}}

\ncm{\CC}{{\bf C}}
\ncm{\ZZ}{{\bf Z}}
\ncm{\BB}{{\bf B}}
\ncm{\LL}{{\bf L}}
\ncm{\JJ}{{\bf J}}
\ncm{\one}{{\bf 1}}

\ncm{\HH}{{H\o H}}
\ncm{\NN}{{\bf N}}
\ncm{\MM}{{\bf M}}
\ncm{\DD}{{\bf D}}
\ncm{\FF}{{\bf F}}
\ncm{\FH}{\FF_{H\o H}}
\ncm{\EE}{{\bf E}}
\ncm{\RR}{{\bf R}}

\def\lcros{\raise1.5pt\hbox{$\scriptstyle\triangleright$}\!
           \raise1.9pt\hbox{$\scriptscriptstyle < \,$}}
\def\cros{\raise1.9pt\hbox{$\scriptscriptstyle  > $}\!
          \raise1.5pt\hbox{$\scriptstyle\triangleleft\,$}}
\def\lef{{\,\hbox{$\textstyle\triangleright$}\,}}
\def\rig{{\,\hbox{$\textstyle\triangleleft$}\,}}
\def\Rep{\mbox{Rep}\,}

\def\Mon{\mbox{\sf Mon}_{R}\,}
\def\mon{\mbox{\sf mon}_{R}\,}
\def\monop{\mbox{\sf mon}_{R_{op}}\,}

\def\Ad{\mbox{Ad}\,}
\def\bra{\langle}
\def\ket{\rangle}
\ncm{\bsn}{\bigskip\noindent}
\renewcommand{\theequation}{\mbox{\arabic{section}.\arabic{equation}}}
\def\o{\otimes}
\def\e{\varepsilon}
\def\x{\times}

\def\be{\beta}
\ncm{\1}{_{(1)}}
\ncm{\2}{_{(2)}}
\ncm{\3}{_{(3)}}
\ncm{\4}{_{(4)}}
\ncm{\note}{\marginpar{\rule[-3mm]{1mm}{6mm}}}
\ncm{\lb}{\label}

\begin{document}

\LARGE

\title{\bf On the Structure of Monodromy Algebras\\[.15cm]
  and Drinfeld Doubles\\[1cm]}
\author{{\sc Florian Nill}%
\footnotemark[1]\\
Institut f\"ur Theoretische Physik,
\\
FU-Berlin, Arnimallee 14, D-14195 Berlin}

\date{Aug. 23, 1996\\ Revised Oct. 6, 1996}
\rncm{\thefootnote}{\fnsymbol{footnote}}
\footnotetext[1]{E-mail: NILL@physik.fu-berlin.de\\
   Supported by the DFG, SFB 288 "Differentialgeometrie und
   Quantenphysik"}
\rncm{\thefootnote}{\arabic{footnote}}
\maketitle
\normalsize
\vspace{-10.5cm}
\begin{flushright}
 q-alg/9609020
\\
To appear in Rev. Math. Phys.
\end{flushright}
\vspace{8.5cm}

\begin{abstract}
We give a review and some new relations on the structure of the
monodromy algebra (also called loop algebra) associated with a
quasitriangular Hopf algebra $H$. It is shown that as an algebra it coincides
with the so-called braided group constructed by S. Majid on the dual
of $H$. Gauge transformations
act on monodromy algebras via the coadjoint action. Applying a
result of Majid, the resulting crossed product is isomorphic
to the Drinfeld double $\D(H)$. Hence, under the so-called factorizability
condition given by N. Reshetikhin and M. Semenov-Tian-Shansky,
both algebras are isomorphic to the algbraic tensor product
$H\o H$. It is indicated that in this way the results of
Alekseev et al. on lattice current algebras are consistent
with the theory of more general Hopf spin chains given by
K. Szlach\'anyi and the author.
In the Appendix the multi-loop algebras $\L_m$ of Alekseev and
Schomerus [AS] are identified with braided tensor products of
monodromy algebras in the sense of Majid, which leads to an
explanation of the ``bosonization formula'' of [AS]
representing $\L_m$ as $H\o\dots\o H$.

\end{abstract}

\sec{Introduction}

Monodromy algebras and Drinfeld doubles have appeared as
quantum symmetry algebras in several models of low dimensional
quantum field theory during the last few years.

Monodromy or loop algebras associated with quasitriangular Hopf
algebras $H$ play an important r\^ole in lattice approaches to
Chern-Simons theory [AGS], topological quantum field theory [AS] and current
algebras on a circle [AFFS].
In these models they commonly
appear as algebras generated by the matrix
elements of quantum holonomy operators $\MM = (M^I_{ab})$
around closed loops.
Their center is spanned by generators obeying the Verlinde algebra
and commuting with all other link operators.
Their represenation theory has
been used in [AGS,AS]
to study quantum Chern-Simons
algebras (``graph algebras")
on Riemann surfaces with coloured punctures.
As a particular input result the authors show that under the
assumption on $H$ to be modular Hopf algebra the representation
category of the associated monodromy algebra coincides with that of
$H$ itself [AS].

In [AFFS] these results have been used and further developped
to study current algebras on a periodic lattice chain. The authors introduce
``gauged loop algebras" $\K$ given as crossed products of a monodromy
algebra with a copy of the gauge quantum group sitting on the
initial ($\equiv$ end) point of the loop. It turns out [AFFS] that the
center of $\K$ already coincides with the center of the whole
lattice current algebra $\A$, thereby governing again its
representation theory.
Under the same assumptions as above it is proven in [AFFS] that
the irreducible representations of $\K$ - and therefore of $\A$ - are
in one-to-one correspondence with those of $H\o H$.
\footnote{This holds true on finite lattices. The discussion
about the infinite lattice limit in [AFFS] is heuristic and
inconclusive in various respects.}

\bigskip
On the other hand, in non-abelian spin systems [SzV] and their
generalizations to Hopf spin sytems
[NSz] it is the Drinfeld double $\D(H)$ [Dr] which plays a
similarly central r\^ole. In fact, it appears as the appropriate
generalization of the familiar order-disorder symmetry group
$G\x G$ known from abelian $G$-spin models. As a
rather complete classification result it has been proven in
[NSz] that in the infinite lattice limit
the category of DHR-superselection sectors of these
models is precisely given by $\Rep(\D(H))$.
This has been done
by constructing an equivalence class of localized transportable coactions
$\rho:\A\to\A\o\D(H)$,
where here $\A$ denotes the quasi-local $C^*$-algebra
generated by all local
observables on the infinite lattice chain.
One may then show that these coactions are universal in
the sense that any charged represenation $\pi_r$ of $\A$ satisfying the
DHR-selection criterion relative to some Haag-dual vacuum
representation $\pi_0$ is equivalent
to a representation $(\pi_0\o D_r)\circ\rho$ where $D_r$ is a
representation of $\D(H)$
and where the correspondence $\pi_r\leftrightarrow
D_r$ is one-to-one on equivalence classes. The associated
statistics operators are $\CC$-matrix valued and are precisely
given by the representation matrices of the
canonical quasi-triangular $R$-matrix of $\D(H)$.
The above coaction is in a way dual to the symmetry action of
$\D(H)$, which leaves the observables invariant and only acts
on the order-disorder fields.

All this confirms with the fact that the
Drinfeld double has also been discovered in 2-dimensional
continuum models to describe the quantum symmetry associated
with nonabelian order and disorder (kink) fields [DPR,M\"u].

\bigskip
Not obviously being related, Majid has
developed his theory of braided groups during the last 5 years
[M3-M7] by also constructing new
algebra and coalgebra structures from a given quasitriangular
Hopf algebra $H$. Particular examples of physical interest are
among others the Sklynanin Algebra and the quantum Lorentz
group [M5]. As a guiding construction principle, braided groups
always appear as (co-)algebras in a category of $H$-modules,
i.e. all structural maps are built so as to become $H$-module
morphisms. In the cases of interest for us one takes $H$ (or
its dual $\hat H$)
itself as the underlying space and considers the (co-)adjoint
H-action. The given structures of $H$ and $\hat H$ must then be
deformed in order to provide $H$-module morphisms under this
action. Based on this philosophy one is also naturally lead to
depart from ordinary tensor products to braided tensor products
[M4]. Majid also provides what he calls ``bosonization
formulas" [M4,M5,M7,M8] by showing that various crossed product
constructions by nontrivial Hopf actions lead to algebras which
are nevertheless isomorphic to the ordinary tensor product of
its constituents. This might already give a first hint to some
hidden relations with the above mentioned results of [AFFS] on
the representation theory of their algebra $\K$.

\bigskip
Digging further into the literature one discovers that already
in 1988 Reshetikhin and Semenov-Tian-Shansky [RS1] have provided
a somewhat related ``bosonization formula" by noting that
under a what they called ``factorizability condition" on the
$R$-matrix of a quasi-triangular Hopf algebra $H$ the
associated Drinfeld double $\D(H)$ is isomorphic to $H\o H$
\footnote{Actually, this is more then just a bosonization, since the copy of
the dual algebra $\hat H\subset \D(H)$ is traded for
a copy of $H$. We will make the relations more explicit by
dividing this transformation
into two steps in Section 7, see also Chapter 7 of [M8]}.
Around the same time Majid had developed his theory of the Drinfeld double
as a double crossed product of $H$ with $\widehat{H_{op}}$. He then
discovered that without any additional conditions
the Drinfeld double of a quasi-triangular Hopf
algebra $H$ is always isomorphic to an ordinary crossed product
$\hat H_R\cros H$
[M2] by constructing a Hopf module action of $H$ on a newly
invented algebra $\hat H_R$, which only latter became to be
described as the ``braided group function algebra" [M4] dual to
the braided group introduced in [M3].

\bigskip
Putting all these pieces together one could already arrive at the conjecture
that Majid`s crossed product $\hat H_R\cros H$
might be the same as the algebra
$\K$ in [AFFS], the latter also being a crossed product
(i.e. $\{\mbox{monodromy algebra}\}\cros H$), and that both algebras are
therefore isomorphic to the Drinfeld double $\D(H)$.
The identification with the representation
category of $H\o H$ would then be a consequence of the results
of [RS1], since
their factorizability condition is in fact nothing but
the invertibility property of the linear ``monodromy map"
$$
\mon: \hat H\ni \varphi \mapsto (\varphi\o id)(R_{op} R)\in H
$$
which was proven (although somewhat hiddenly) to hold for
modular Hopf algebras in Section 5.2 of [AS].

\bigskip
This picture may be further supported by the observation, that
apart from the periodic boundary condition the lattice current
algebra of [AFFS] is in fact a reformulation of the Hopf spin
chain of [NSz].
To see this let me shortly review that the latter is defined by
placing a copy of a finite dimensional
Hopf algebra $\A_{2i}\cong H$ on each site and a copy of the dual algebra
$\A_{2i+1}\cong\hat H$ on each link of a one dimensional lattice.
Non-vanishing commutation relations are then postulated only on neighbouring
site-link pairs, where one requires [NSz]
\bea
A_{2i} (a)A_{2i-1}(\varphi) &=&
A_{2i-1}(\varphi\1)\,\bra\varphi\2\mid a\1\ket\,A_{2i}
(a\2)\label{1.1}\\
A_{2i+1}(\varphi)A_{2i}(a)&=&A_{2i} (a\1)\,\bra a\2\mid\varphi\1\ket\,A_{2i+1}
(\varphi\2)\label{1.2}
\eea
Here sites (links) are numbered by even (odd)
integers and  $H\ni a\mapsto A_{2i}(a)\in \A_{2i}\subset\A,\
\hat H\ni\varphi\mapsto A_{2i+1}(a)\in \A_{2i+1}\subset\
A$ denote the
embeddings of the single site (link) algebras into the global quantum
chain $\A$. As usual,
$a\1\o a\2 \equiv \sum_i a\1^i \o a\2^i=\Delta(a)$ denotes
the coproduct on $H$ (and analogously on $\hat H$) and
$\bra a\mid\varphi\ket\equiv
\bra\varphi\mid a\ket$ denotes the dual pairing $H\o \hat H\to \CC$.

\smallskip
Let now $R\in H\o H$ be quasitriangular and put
$\RR_{2i} := (id\o A_{2i})(R_{op}^{-1})
\in H\o \A_{2i}$, where $R_{op}$ denotes the image of $R$ under
the permutation of tensor factors. Let us also introduce
\beq \label{1.3}
\LL_{2i+1} := (id\o A_{2i+1})(\EE) \in H\o \A_{2i+1}
\eeq
where $\EE =\sum e_\nu\o e^\nu\in H\o\hat H$
is the canonical element given by
any basis $e_\nu \in H$ with dual basis $e^\nu \in\hat H$.
One may then define ``lattice currents"
\beq\label{1.4}
\JJ_{2i+1} := \RR_{2i}^{-1} \LL_{2i+1} \in H\otimes \A
\eeq
which are immediately verified to satisfy the lattice current algebra of
[AFFS]
\bea
[\one_H\o A_{2i}(a)] \JJ_{2i-1} &=& \JJ_{2i-1} [a\1\o A_{2i} (a\2)]
\label{1.5}\\
\JJ_{2i+1} [\one_H\o A_{2i} (a)] &=& [a\1 \otimes A_{2i} (a\2)] \JJ_{2i+1}
\label{1.6}\\
\JJ^{13}_{2i+1} \JJ^{23}_{2i+1}
&=& R^{12} (\Delta\o id)(\JJ_{2i+1}) \label{1.7}\\
\JJ^{13}_{2i-1}R^{12} \JJ^{23}_{2i+1} &=&
\JJ^{23}_{2i+1} \JJ^{13}_{2i-1}\label{1.8}
\eea
where the last two lines are
understood as identities in $H\o H\o \A$, the upper
 indices indicating the canonical embeddings of tensor factors.

Hence, under the additional data given by a quasitriangular $R$ the lattice
algebras of [NSz] and [AFFS] are isomorphic. Moreover, periodic boundary
conditions may be imposed in both versions by identifying
$\A_0\equiv \A_{2N}$.

\bigskip
The discovery of this isomorphism was the starting motivation
for the present work. In fact, it seemed hard to believe that the
different results on the representation theory of this model
(i.e. isomorphic to $\Rep\K\equiv\Rep H\o H$ in [AFFS] and to
$\Rep\D(H)$ in [NSz]) should be caused by the different
approaches of the authors and/or the different boundary conditions.

As it turns out now, the true explanation is indeed given by the
chain of isomorphisms mentioned above
$$
\K\cong\hat H_R\cros H\cong\D(H)\cong H\o H
$$
(the last one holding provided $\mon$ is invertible).
Moreover, all of this could in principle
be deduced from the literature as quoted,
provided one realizes that the braided group function
algebra $\hat H_R$ of Majid is in fact the same object as the
monodromy algebra introduced in [AGS,AS,AFFS].
So with this point of view the paper could end with just proving this last
statement.

\bigskip
However I would say that even the experts agree
that the relevant informations are spread in an rather
unrelated way over the
literature and often have very different appearances and
notations making it hard to compare or even identify them.
Due to the importance of these structures in the
various mathematical and physical contexts mentioned above I therefore
consider it worthwhile to give a selfcontained  and
unifying account of the whole story here.

So let me
start in Section 2 by identifying the monodromy algebras of [AGS,AS,AFFS]
as (duals of) braided groups $\hat H_R$ in the sense of Majid [M4,M5].
In Section
3 I will then review Majid's construction of a coadjoint $H$-action on
$\hat H_R$ and identify the resulting crossed product $\M_R(H):=\hat
H_R\cros H$
with the gauged monodromy algebra $\K$ given in Definition 5 of
[AFFS]. Section 4  gives a complete proof of an observation
of [AFFS] on the existence of commuting left and right
monodromies inside $\M_R(H)$. In the course of this proof we find
a copy of $\hat H$ naturally embedded as a subalgebra in
$\M_R(H)$. In this way we recover as a Corollary in Section 5
the isomorphism $\M_R(H)\cong\D(H)$
previously obtained in [M2]. In Section 6 we use methods
of [M5] to show that the monodromy map $\mon$ extends to a
homomorphism of crossed products $\Mon : \M_R(H)\to H\cros_{Ad}H$
which becomes an isomorphism iff $\mon$ is bijective (here $Ad$
denotes the adjoint action of $H$ on itself).
Since by a simple bosonization formula the crossed product
$H\cros_{Ad}H$ is naturally isomorphic to the algebraic tensor
product $H\o_{alg} H$, this clarifies the results of [AFFS] on the
representation theory of $\K\equiv\M_R(H)$ without relying on
any semisimplicity assumptions.

In Section 7 we proceed to the Hopf algebra
structure induced on
$\M_R(H)\equiv \K$ by the isomorphism with $\D(H)$ and investigate its
image on $H\o_{alg}H\cong H\cros_{Ad}H$ under the map $\Mon$. We
define on $H\o_{alg}H$ a cocycle deformed
coproduct $\delta_\HH$ and an associated quasitriangular
R-matrix $\R$. In this way the map $\D(H)\cong\M_R(H)\to H\o_{alg}H$
induced by $\Mon$ provides a homomorphism (isomorphism, if
$\mon$ is bijective) of quasitriangular Hopf algebras, which turns out to
coincide with the one already given in [RS1].
It is conjectured that after some corrections
\footnote{There were some inconsistencies in [AFFS], see Section 7}
the Hopf algebra structure obtained on
$\K$ by [AFFS]
should agree with these results.
Furthermore, with these corrections
the $\K$-coaction on $\A$ discovered by [AFFS]
should also be expected to become an example from the class of universal
localized $\D(H)$-coactions on $\A$ established in [NSz].

A more detailed account of this will be given elsewhere, where
it will also be shown
that, independently of the existence of $R$,
the model (1.1), (1.2) with periodic
boundary conditions always contains $\D(H)$ as a global subalgebra
governing its representation theory, similarly as in [AFFS].

The Appendix concludes with an independent investigation of
the multi-loop algebras $\L_m$ of [AS]. They are shown to be
braided tensor products in the sense of
Majid of one-loop ($\equiv$ monodromy) algebras.
Hence, using bosonization ideas of [M5-M8], they are
naturally homomorphic (isomorphic, if $\mon$ is bijective) to
the algebraic tensor product $H^{\o m}$. This explains the
representation theory of $\L_m$ given in [AS].

\bigskip
In summary, I would once more like to emphasize that in most parts this
paper should be considered as a review, unifying different
notations and approaches and treating seemingly unrelated
results within a common formalism.
In particular, Majid's notion of braided groups is brought
together with the techniques of generating matrices as advocated
by the St. Petersburg school.
In view of the fact that the various relations between
monodromy algebras, Drinfeld doubles, braided groups and
bosonization formulas seemingly have been overlooked in
[AGS,AS,AFFS], this article is hoped to add a piece of clarity
and also to serve a broader community working in the field.

\bigskip
I should point out that a rather extensive
review on braided groups and braided
algebras can be found in [M7] and the relation with the
Drinfeld double is also reviewed in [M8], however without
relating it to monodromy algebras and without using the
techniques of generating matrices.

\bigskip
Throughout this paper all algebras are taken to be finite
dimensional, but no assumptions on semi-simplicity nor on
$*$-structures are made. The reader is supposed to be familiar
with standard Hopf algebra theory and notation, see e.g. [Sw].
Elements of $H$ will be denoted by $a,b,c,\dots$ and elements
of $\hat H$ by $\varphi,\psi,\chi,\dots$. The counit is denoted
by $\e :H\to\CC$ and the antipode by $S:H\to H$.

\bsn
{\bf Acknowledgements:}
I would like to thank A. Alekseev and V. Schomerus for
discussions on their results in [AGS,AS,AFFS]. I am also
indepted to S. Majid for helpful comments and for bringing
refs. [M6-M8] to my attention.

\sec{Monodromy Algebras and Braided Groups}

Defining relations of so-called monodromy or loop algebras have
been given in [AGS,AS], where they appear as algebras generated
by the matrix elements of quantum holonomy operators around
closed loops. For an earlier version based on quadratic
relations see also [RS2]. In order to establish that they are indeed
well defined associative algebras we start by
identifying their dual coalgebras with the
braided groups introduced by Majid in [M3].

\bsn
{\bf Proposition 2.1} [M3]
{\sl Let $R=\sum_{i} x^i\o y^i\in H\o H$ be quasitriangular
with respect to $\Delta$ and let $\Delta_R:H\to H\o H$ be a
deformed coproduct given by
\beq \label{2.1}
\Delta_R(a)= \sum_{i,j} (x^i\o S(y^j)y^i) \Delta(a)(x^j\o\one),\ a\in H
\eeq
Then $\Delta_R$ is coassociative with counit $\e$.}

\bsn
{\bf Proof:} Recall that quasitriangularity of $R$ means [Dr2]
\bea
(\Delta\o id)(R) &=& R^{13} R^{23} \label{2.2} \\
(id\o\Delta)(R) &=& R^{13} R^{12} \label{2.3} \\
R\,\Delta(a)\,R^{-1} &=& \Delta_{op}(a) \label{2.4}
\eea
This implies $(\e\o id)(R)=(id\o\e)(R)=\one$ and therefore $\e$ is a counit for
$\Delta_R$. Moreover, (2.2)-(2.4) also imply the {\em cocycle
  property} [Dr2]
\beq\label{2.5}
R^{12} (\Delta\o id) (R) = R^{23} (id\o\Delta)(R)
\eeq
which is sufficient for $\Delta':H\to H\o H$
\beq\label {2.6}
\Delta'(a):=R\Delta(a)
\eeq
to be coassociative.  Now $\Delta'$ obeys
\beq\label{2.7}
\Delta'(ab) = \Delta' (a)\Delta(b) =\Delta_{op} (a) \Delta'(b)
\eeq
for all $a,b\in H$, from which Proposition 2.1 follows by using

\bsn
{\bf Lemma 2.2:}
{\sl Let $\Delta':H\to H\o H$ be a coassociative coproduct satisfying
(2.7). Then
$
\Delta_R(a) := \sum_j (\one \otimes S(y^j)) \Delta' (a)(x^j\otimes \one)
$
is also coassociative.
}

\bsn
{\bf Proof:} Using (2.7) and (2.2) we compute for $a\in H$
\beanon
(\Delta_R\o id)(\Delta_R(a))
&=& \sum_{i,j} [\one\o S(y^i)\o S(y^j)] \Delta^{'(2)} (a)[x\1^j x^i\o x\2^j
\o\one]\\
&=& \sum_{i,j,k} [\one\o S(y^i)\o S(y^j y^k)] \Delta^{'(2)}(a)
[x^j x^i\o x^k\o \one]
\eeanon
where $\Delta^{'(2)} = (\Delta'\o id)\circ \Delta'
=(id\o\Delta')\circ\Delta'$.
Similarly, using (2.7) and (2.3),
\beanon
(id\o\Delta_R)(\Delta_R(a))
&=& \sum_{i,j} [\one \o S(y\1^i)\o S(y^k)S(y\2^i)]
\Delta^{'(2)} (a) [x^i\o x^k\o \one]\\
&=& \sum_{i,j,k} [\one\o S(y^i)\o S(y^k)S(y^j)] \Delta^{'(2)} (a)
[x^jx^i \o x^k \o \one]\\
&=& (\Delta_R\o id)\Delta_R(a)\qquad\qed
\eeanon
 It has been remarked in equ. (32) of [M6] that the passage
  from $\Delta'$ to $\Delta_R$ can also be viewed as a cocycle
  transformation in a generalized sense, i.e. by considering
according to (2.7) $(H,\Delta')$ as an $(H\o H)$-Hopf right-module coalgebra
 with right action $a\rig (b\o c):=S(b)ac$, and by acting with
the $H\o H$-right cocycle $R_{23}\in (H\o H)\o (H\o H)$ (see
Lemma 7.1) to transform $\Delta'\to\Delta_R$.

\bsn
Up to a change of conventions (i.e. replacing $(H,\Delta,S,R)$
by $(H,\Delta_{op},S^{-1},R_{op})$\,)
the deformed coproduct \no{2.1} has first appeared in [M3],
providing $H$ with the structure of what the author calls
a ``braided group''.
One should also point out here that the deformed coproducts $\Delta'$ and
$\Delta_R$ are no longer algebra homomorphisms. However, as has been
emphasized by Majid [M3,M4], there exists a natural deformed algebra
structure on $H\o H$, the ``braided tensor product'' $H\o_R H$ (see
Appendix A), such that $\Delta_R :H_{cop}\to H_{cop}\o_{R_{op}}
H_{cop}$ does provide an algebra homomorphism (see equ.\no{A.11}), where
$H_{cop}$ is the Hopf algebra $H$ with opposite coproduct
$\Delta_{op}$.

Following Theorem 4.1 of [M4], see also equ. (8) of [M5],
we now pass to the algebra dual to $(H,\Delta_R)$.

\bsn
{\bf Definition 2.3} [M4]
{\sl The ``braided group function algebra'' $\hat H_R$ is defined to be the
vector space $\hat H$ with multiplication $\x_R : \hat H\o\hat H\to\hat H$
induced by $\Delta_R$, i.e.
\beq \label{2.8}
\bra \varphi  \x_R \psi\mid a\ket :=\bra\varphi\o\psi\mid\Delta_R(a)\ket
\eeq
for $\varphi,\psi\in\hat H$ and $a\in H$.}

\bsn
Note that $\Delta_R$ being coassociative, the new
product $\x_R$ on
$\hat H$ is clearly associative with unit given by $\e \in \hat H$.
Let me now show that the Definition 2.3 coincides with the
description of monodromy algebras given by
[AGS,AS,AFFS] in terms of a generating matrix $\MM$.

\bigskip
In what follows, $\A$ always denotes an arbitrary target algebra.
The basic idea behind generating matrices is given by the observation that if
$\hat\B$ is the algebra dual to some given finite dimensional coalgebra
$(\B,\Delta_\B,\e_\B)$, then the relation
$f(\varphi)=(\varphi\o id)(\FF_\A),\ \varphi\in\hat\BB$,
provides a one-to-one correspondence between algebra homomorphisms
$f:\hat\B\to\A$ and ``generating matrices" $\FF_\A\in\B\otimes \A$ obeying
\beq\label{2.9}
\FF_\A^{13} \FF_\A^{23} = (\Delta_\B \o id)(\FF_\A)\in \B\o\B\o\A
\eeq
Clearly, $f$ is unital iff $\FF_\A$ is unital in
the sense $(\e_\B\o id)(\FF_\A) = \one_\A$.
If $\B$, and therefore $\hat\B$, are bialgebras, then $\B$ has an
antipode $S_\B$ if and only if all unital generating matrices
are invertible, the relation being given by
\beq\label{2.9b}
\FF_\A^{-1} = (S_\B\o id_\A)(\FF_\A)
\eeq
If also $\A$ is a bialgebra then $f$ defines a bialgebra
homomorphism if and only if in addition to (2.9) $\FF_\A$ satisfies
\beq\label{2.10}
\FF_\A^{12} \FF_\A^{13} =(id\o\Delta_\A)(\FF_\A)\in\B\o\A\o\A
\eeq
Putting $\A=\hat\B$ and $f=id_{\hat\B}$ one obtains $\FF \equiv \EE=\sum_\nu
e_\nu\o e^\nu \in \B\o \hat\B$,
which justifies the statement that the algebra
$\hat\B$ is generated by the entries $(\varphi\o id)(\FF),\
\varphi\in\hat\B$,
with abstract relations (2.9).
Applying all this to the algebra $\hat H_R$ we get

\bsn
{\bf Proposition 2.4}
{\sl The relation $f(\varphi)=(\varphi\o id)(\MM_\A)$ provides
a one-to-one
correspondence between algebra homomorphisms $f:\hat H_R\to \A$ and
elements $\MM_\A\in H\o \A$ satisfying}
\beq \label{2.11}
\MM_\A^{13} R^{12} \MM_\A^{23} = R^{12} (\Delta\o id)(\MM_\A)\in H\otimes H
\otimes\A
\eeq
{\bf Proof:} According to (\ref{2.8})
and (2.9) we have to verify that \no{2.11} is
equivalent to
\beq \label{2.12}
\MM_\A^{13} \MM_\A^{23} = (\Delta_R\o id)(\MM_\A)
\eeq
To this end we use (2.3) to get the identity
\beq\label{2.14}
\sum_{i,j} x^ix^j\o S(y^j)y^i
=\sum_{i} x^i\o S(y\1^i) y\2^i = \one\o\one
\eeq
Hence \no{2.11} implies
$$
\MM^{13}_\A \MM^{23}_\A =\sum_j [\one\o S(y^j)\o\one]\MM_\A^{13} R^{12}
\MM_\A^{23} [x^j\o\one \o\one]
=(\Delta_R\o id)(\MM_\A)
$$
Conversely, \no{2.12} implies
$$
\MM_\A^{13} R^{12} \MM_\A^{23} = \sum_i [\one\o y^j\o\one] (\Delta_R\o id)
(\MM_\A)[x^j\o\one\o\one]
= R^{12} (\Delta\o id)(\MM_\A)
$$
where, similarly as above, we have used
\beq\label{2.15}
\sum_{i,j} x^ix^j\o y^j S(x^i) =
\sum_i x^i\o y\1^i S(y\2^i) =\one\o\one
\eeq
\qed

\bsn
Putting $\A = \hat H_R$ and $f=id$,
equation (\ref{2.11}) coincides with the defining ``monodromy
relations'' given in [AFFS]. The convention of [AGS,AS] differs by a
multiplication of $\MM_\A$ with $(\kappa^{-1}\o\one_\A)$, where
$\kappa\in H$ is a certain central square root of a ribbon element
associated with $R$.

Hence, from now on we will call $\x_R$ the {\em monodromy product} on
$\hat H$ and $\hat H_R$ the {\em monodromy algebra} associated with $(H,R)$.

\sec{Gauged Monodromy Algebras}

Following Majid [M4,M5]
we now provide a coadjoint action of $H$ on $\hat H_R$ which becomes
a Hopf module action with respect to the monodromy product $\x_R$. Again we
start with the dual point of view by first recalling that
$\Delta_R$ is an $H$-module map with respect to the right
adjoint action of $H$ on itself (this is a basic ingredient
in Majid's construction of braided groups by what he calls
``transmutation'').

\bsn
{\bf Lemma 3.1}
{\sl For $a,b\in H$ we have}
\beq\label{3.1}
\Delta_R(S(a\1)ba\2)=[S(a\1)\o S(a\3)] \Delta_R(b)[a\2\o a\4]
\eeq
{\bf Proof:} By straight forward verification, using (2.4), (2.7) and the
identity $\Delta\circ S=(S\o S)\circ \Delta_{op}$.

\bsn
{\bf Corollary 3.2}
{\sl Let $\lef :H\o \hat H_R\ni a\o\varphi \mapsto a\lef\varphi
\in \hat H_R$ be given by
\beq\label{3.2}
\bra a\lef\varphi\mid b\ket := \bra\varphi\mid S(a\1)ba\2\ket
\eeq
where $a,b\in H$ and $\varphi\in\hat H_R$. Then $\lef$ defines a Hopf module
left action of $H$ on $\hat H_R$.}

\bsn
{\bf Proof:} The identities $a\lef (b\lef \varphi)=(ab)\lef\varphi,
\ \one \lef\varphi=\varphi$ and $a\lef\e =\e(a)\e$ are obvious.
The Hopf module property
\beq\label{3.3}
a\lef (\varphi\x_R\psi) =(a\1\lef\varphi)\x_R(a\2\lef\psi)
\eeq
follows by pairing both sides with $b\in H$ and using (\ref{2.8}) and (3.1).
\qed

\bsn
{\bf Definition 3.3}
{\sl The ``gauged monodromy algebra'' $\M_R(H)$ is defined
to be the crossed product}
$$
\M_R(H):= \hat H_R\cros H
$$

\bsn
Thus, recalling the definition of a crossed product, $\M_R(H)$ is the linear
space $\hat H\o H$ with associative multiplication
\beq\label{3.4}
(\varphi\o a)(\psi\o b) := (\varphi \x_R(a\1 \lef\psi)\o a\2 b),
\eeq
where $\varphi,\psi\in\hat H$ and $a,b\in H$.
To distinguish the new algebraic structure we denote the images of $\hat H_R$
and $H$ in $\M_R(H)$ by $i_M(a)\equiv \e\o a,\ a\in H$,
and $M(\varphi)\equiv\varphi\o \one_H,\ \varphi\in\hat H$.
We also denote $\MM:=\sum e_\nu\o M(e^\nu)\in H\o \M_R(H)$.

The crossed product $\M_R(H)$ has first appeared in [M2,M5].
Using the generating-matrix formalism
we now show that $\M_R(H)$ coincides with the
algebra $\K$ given in Definition 5 of [AFFS] (see also [AGS,AS])
as the subalgebra generated inside a current algebra by a monodromy
operator $\MM$ together with the algebra
of gauge transformations sitting at its initial ($\equiv$ end) point.

\bsn
{\bf Proposition 3.4}
{\sl Let $f:H\to\A$ be an algebra homomorphism. Then the relation
$$
f_M(\varphi\o a):= (\varphi\o id)(\MM_\A) f(a)
$$
provides a one-to-one
correspondence between algebra
homomorphisms $f_M:\M_R(H)\to\A$ extending $f$ and
elements $\MM_\A \in H\o \A$ satisfying \no{2.11} together with
$(\varepsilon\o id)(\MM_\A)= f(\one_H)$ and}
\beq\label{3.5}
[a\1\o f(a\2)] \MM_\A=\MM_\A [a\1\o f(a\2)],~~\forall a\in H.
\eeq
{\bf Proof:} %
Since $\one_{\M_R(H)}=\varepsilon\o\one_H$, the condition
$(\varepsilon\o id)(\MM_\A)= f(\one_H)$ is equivalent to
$f_M(\one_{\M_R(H)})=f(\one_H)$.
We are left to show that (3.5) holds if and only if $f_M$ respects
(3.4), or equivalently, if and only if
\beq\label{3.6}
f(a\1)(\varphi \o id_\A)(\MM_\A) f(S(a\2)) =((a\lef\varphi)\o id)(\MM_\A)
\eeq
for all $\varphi \in\hat H$ and all $a\in H$. Now (3.6) is equivalent to
\beq\label{3.7}
(\one_H\o f(a\1))\MM_\A (\one_H\o f(S(a\2))=(S(a\1)\o \one_\A)(\MM_\A)
(a\2 \o\one_\A)
\eeq
which further implies
\beanon
[a\1\o f(a\2)] \M_\A &=& [a\1\o f(a\2)] \MM_\A [\one_H\o f(S(a\3)a\4)]\\
&=& [a\1 S(a\2)\o \one_\A] \MM_\A [a\3 \o f(a\4)]\\
&=& \MM_\A[a\2 \o f(a\2)]
\eeanon
and therefore (\ref{3.5}). Conversely, given (\ref{3.5}) we conclude
\beanon
[\one_H\o f(a\1)] \MM_\A[\one_H\o f(S(a\2))] &=&
 [S(a\1)a\2\o f(a\3)] \MM_\A[\one_H\o f(S(a\4))]\\
&=& [S(a\1)\o\one_\A]\MM_\A[a\2 \o f(a\3 S(a\4))]\\
&=& [S(a\1)\o\one_\A] \MM_\A[a\2 \o\one_\A]
\eeanon
and therefore (3.7). \qed

\sec{Left and Right Monodromies}

Here we review the observation of [AFFS] on a right monodromy algebra
$\hat H_R^r \subset \M_R(H)$
commuting with $\hat H_R^\ell \equiv\hat H_R\subset
\M_R(H)$.
First we define $\hat H_R^r$ abstractly as in Definition 2.3, where we replace
$(H,\Delta,R,S)$ by $(H,\Delta_{op}, R_{op}, S^{-1})$, i.e. the
quasitriangular  Hopf
algebra $H$ with opposite coproduct.
Similarly, to define a left Hopf module action of
$H_{cop}\equiv (H,\Delta_{op})$ on
$\hat H_R^r$ we replace (3.2) by
\beq\label{4.1}
\bra a\lef'\varphi\mid b\ket := \bra\varphi \mid S^{-1} (a\2)ba\1\ket
\eeq
leading to a crossed product
\beq\label{4.2}
\M_R^r(H):=\hat H_R^r \cros H_{cop}
\eeq
for which Propositions 2.4 and 3.4 hold analogously,
with $(\Delta,R)$ replaced by $(\Delta_{op}, R_{op})$.
Putting in particular $\A=\M_R^\ell (H)\equiv
\M_R(H)$ we  now get an isomorphism $\M_R^r(H)\to \M_R^\ell (H)$ restricting
to the identity on $H$ by choosing in Proposition 3.4 $\MM_\A = \MM^r$,
\beq\label{4.3}
\MM^r:= \RR\,(\MM^\ell)^{-1} \RR_{op} \in H\o\M_R^\ell(H)
\eeq
where $\MM^\ell \equiv
\MM:=\sum e_\nu\o M(e^\nu)\in H\o\M^\ell_R(H)$ and where
$\RR_{(op)}:=(id\o i_M)(R_{(op)})\in H\o\M^\ell_R(H)$.
This is the formula of [AFFS], which is proven by equ.
(\ref{4.6}) below. Note however that for (4.3) to be
well defined we first have to assure that $\MM^\ell$
is invertible in $H\o \M_R^\ell(H)$.

\bsn
{\bf Proposition 4.1}
i) {\sl For any target algebra $\A$ let $f:H\to\A$ be a unital
homomorphism and let $\MM_\A^\ell\equiv\MM_\A\in H\otimes \A$
be unital and satisfy \no{2.11} and (3.5). Then}
\beq\label{4.4}
\MM_\A^{-1} = (S\o id)\left((id\o f)(R_{op}^{-1})\MM_\A\right)
(id\o f)(R_{op}^{-1})
\in H\o\A
\eeq
ii) {\sl Putting
$\MM_\A^r := (id\o f)(R) (\MM_\A^\ell)^{-1} (id\o f) (R_{op})$
we have}
\bea
[a\2\o f(a\1)] \MM_\A^r &=& \MM_\A^r [a\2\o f(a\1)]  \label{4.5} \\
(\MM_\A^r)^{13} R^{21} (\MM_\A^r)^{23} &=&
R^{21}(\Delta_{op}\o id)(\MM_\A^r)  \label{4.6} \\
(\MM_\A^\ell)^{13} (\MM_\A^r) ^{23} &=&
(\MM_\A^r)^{23} (\MM_\A^\ell)^{13} \label{4.7}
\eea
{\bf Proof:} i) Putting
\beq\label{4.8}
\DD_\A := (id \o f)(R_{op}^{-1})\MM_\A \in H\otimes \A
\eeq
and using (2.5) we have in $H\o H\o \A$
\beq\label{4.9}
R^{12} (\Delta\o id)(\MM_\A)=[(R\o\one)(\Delta \o id)(R)]^{312} (\Delta\o id)
(\DD_\A)
\eeq
On the other hand (2.2) and (2.4) imply
\bea\label{4.10}
\MM_\A^{13} R^{12} \MM_\A^{23} &=& \MM_\A^{13} (\Delta \o id)(R)^{132}
\DD_\A^{23} \nonumber\\
&=& [(\Delta \o id)(R) (R_{op}\o \one)]^{132} \DD_\A^{13} \DD_\A^{23}\nonumber \\
&=& [(R\o\one)(\Delta\o id)(R)]^{312} \DD_\A^{13} \DD_\A^{23}
\eea
where in the second line we have used (3.5) and where by
a convenient abuse of
notation we have dropped the symbol $f$. Hence comparing (4.9) and
(4.10) and using \no{2.11} we conclude
\beq\label{4.11}
\DD_\A^{13} \DD_\A^{23} =(\Delta\o id_\A) (\DD_\A)
\eeq
Now $\DD_\A$ is unital since $f$ and $\MM_\A$ are unital and, therefore,
\no{4.11} implies by \no{2.9b}
\beq\label{4.12}
\DD_\A^{-1} =(S\o id_\A)(\DD_\A)
\eeq
from which (4.4) follows.

ii) Equ.(4.5) immediately follows from (2.4) and (3.5).
To prove (4.6) we compute,
omitting the symbol $f$,
\beanon
(\MM_\A^r)^{23} R^{12}(\MM_\A^r)^{13} &=& R^{23} (\MM_\A^{-1})^{23}
R^{32} R^{12} R^{13} (\MM_\A^{-1})^{13} R^{31}\\
&=& R^{23}(\M_\A^{-1})^{23} (id\o\Delta)(R) (R^{-1})^{12}
(\Delta\o id)(R)^{132}(\MM_\A^{-1})^{13} R^{31} \\
&=& R^{12} (\Delta \o id)(R) [(\MM_\A^{-1})^{23} (R^{-1})^{12}
(\MM_\A^{-1})^{13}](\Delta\o id)(R)^{132} R^{31}\\
&=& R^{12}(\Delta\o id)(\MM_\A^r R_{op}^{-1}) (R^{12})^{-1}
R^{12} R^{32} R^{31}\\
&=& R^{12} (\Delta\o id)(\MM_\A^r)
\eeanon
Here we have used (2.2)-(2.4) in the second line, (3.5) and
(\ref{2.5}) in the third
line, (2.2), (\ref{4.3}) and the inverse of \no{2.11} in the fourth line and
again (\ref{2.3}) in the last
line. Applying the permutation (12) $\to$ (21) to both sides this proves
(4.6).
Finally, to prove (4.7) we put $\MM\equiv M_\A^\ell$
and $\Omega\equiv (M_\A^r)^{-1}$ and compute
\beanon
\MM^{13} \Omega^{23} &=&
\MM^{13} (R^{12}R^{32})^{-1} R^{12} \MM^{23} (R^{23})^{-1}\\
&=& (\Delta\o id)(R^{-1})^{132} R^{12} (\Delta \o id)(\MM) (R^{23})^{-1}\\
&=& (R^{32})^{-1} (\Delta\o id)(\MM)(R^{23})^{-1}
\eeanon
where in the second line we have used (\ref{2.2}), (\ref{3.5}) and
(\ref{2.11}), and in the third line again (\ref{2.2}).
Similarly one obtains
\bea\label{4.14}
\Omega^{23} \MM^{13} &=& (R ^{32})^{-1} \MM^{23} R^{21} (R^{23} R^{21})^{-1}
\MM^{13}\nonumber\\
&=& (R^{32})^{-1} R^{21} (\Delta_{op} \o id)(\MM)
(id\o \Delta)(R^{-1})^{213}\nonumber\\
&=& (R^{32})^{-1} R^{21} (\Delta_{op} \o id)(\MM)(R^{21})^{-1}
(R^{23})^{-1}\nonumber\\
&=& M^{13} \Omega^{23}
\eea
which proves (4.7). \qed

\bsn
Note that (4.7) implies that the image
of $\hat H_R^r$ in $\M_R^\ell(H)$ indeed
commutes with $\hat H_R^\ell$. Part ii) of Proposition 4.1 has been taken
over from [AFFS].

\sec{The Drinfeld Double}

The Drinfeld double $\D(H)$ (also called quantum double) has been introduced
in [Dr1] and is meanwhile well understood as a
double crossed product of $H$ with
$\hat H$ [M1,M8], see also [RS1,K].

\bsn
{\bf Definition 5.1}
{\sl The Drinfeld double $\D(H)$ over a finite
dimensional Hopf algebra $H$ is the
linear space $\hat H\o H$ with multiplication given for $a,b\in H$ and
$\varphi,\psi\in\hat H$ by}
\beq\label{5.1}
(\varphi\o a)(\psi\o b):=(\varphi\psi\2 \o a\2 b)\
\bra a\1\mid\psi\3 \ket\ \bra\psi\1\mid S^{-1} (a\3)\ket
\eeq

\bsn
Putting $i_D(a):=(\one_{\hat H} \o a)$ and $D(\varphi) :=(\varphi\o \one_H)$
one can rewrite this equivalently as [NSz]
\bea
i_D(a)i_D(b) &=& i_D(ab) \label{5.2} \\
D(\varphi)D(\psi) &=& D(\varphi\psi) \label{5.3} \\
D(\varphi\1)\ \bra\varphi\2\mid a\1\ket\ i_D(a\2) &=& i_D(a\1)\ \bra a\2\mid
\varphi\1\ket\ D(\varphi\2) \label{5.4}
\eea
Hence, as an algebra $\D(H)=\D(\hat H)$.
Based on a more general setting given by [Ra] it has first been
noticed in [M2] (see also Proposition 4.1 of [M5]) that $\D(H)$ is in fact
isomorphic to $\M_R(H)$ for all quasitriangular $R\in H\o H$.

\bsn
Using our formalism of generating matrices let me now
demonstrate that this result reduces to a Corollary of the
calculation leading to \no{4.11}.
First we need an analogue of Proposition 3.4.
Introducing $\DD:=\sum_\nu e_\nu\o D(e^\nu)\in H\o \D(H)$ we note that (5.3)
is equivalent to
$$
\DD^{13} \DD^{23} =(\Delta \o id)(\DD)
$$
and (5.4) is equivalent to
$$
\DD[a\1\o i_D(a\2)] =[a\2\o i_D(a\1)]\DD
$$
More generally this leads to

\bsn
{\bf Lemma 5.2}
{\sl Let $f:H\to \A$ be an algebra
homomorphism. Then the relation
$$
f_D(\varphi\o a) =(\varphi\o id_\A)(\DD_\A) f(a),\ a\in H,\varphi\in\hat H
$$
provides a one-to-one correspondence between algebra homomorphisms
$f_D:\D(H)\to \A$ extending $f$ and elements $\DD_\A\in H\o \A$
obeying $(\varepsilon\o id)(\DD_\A)=f(\one_H)$ and
}
\bea
\DD_\A^{13} \DD_\A^{23} &=& (\Delta\o id)(\DD_\A)\label{5.5}\\
\DD_\A[a\1\o f(a\2)] &=& [a\2\o f(a\1)] \DD_\A,~~a\in H\label{5.6}
\eea

\bsn
{\bf Proof:}
First, by (2.9) the $\DD_\A$'s satisfying (\ref{5.5}) are in one-to-one
correspondence with homomorphisms $\hat H\to \A$.
Since $\one_{\D(H)}=i_D(\one_H)$, the condition $(\varepsilon\o
id)(\DD_\A)=f(\one_H)$ is equivalent to $f_D(\one_{\D(H)})=f(\one_H)$.
We are left to show that
$f_D$ respects (\ref{5.4}) if and only if $\DD_\A$ obeys (\ref{5.6}), which may
immediately be realized by applying $(\varphi\o id_\A)$ to
(\ref{5.6}).\qed

\bsn
Inspired by (4.8) and (4.11) we now put $\A=\M_R(H),\
f=i_M :H\to\M_R(H)$ the
canonical embedding and
\beq\label{5.7}
\DD_M:=\RR_{op}^{-1} \MM\in H\o\M_R(H)
\eeq
where $\RR_{op}:=(id\o i_M)(R_{op})\in H\o\M_R(H)$ and
$\MM=\sum e_\nu\o M(e^\nu)$ as before.
Then we have

\bsn
{\bf Corollary 5.3} [M2,M5]
{\sl The element $\DD_M$ (5.7) defines an algebra isomorphism
$\lambda_R:\D(H)\to\M_R(H)$
restricting to the identity on $H$ by putting}
$$\lambda_R(\varphi\o a):=(\varphi\o id)(\DD_M)\,i_M(a)\ .
$$

\bsn
{\bf Proof:} We apply Lemma 5.2. Equ. (5.5) has already been verified in
(4.11). Equ. (5.6) follows from (3.5) and (2.4). Finally, $\lambda_R$ is
invertible with $\lambda_R^{-1}$ given according to Proposition 3.4 by
\beq \label{5.8}
\MM_D:=(id\o i_D)(R_{op})\DD \in H\o \D(H)
\eeq
\qed

\bsn
Denoting $\MM_D^\ell \equiv \MM_D$ and looking at (4.3) we also get an
immediate formula for a right monodromy $\MM_D^r\in H \o \D(H)$
\beq\label{5.9}
\MM_D^r := (id\o i_D)(R) \DD^{-1}
\eeq
Equ. (4.7) then implies that the subalgebras $\bra \hat H\o id\mid
\MM_D^r\ket \cong \hat H_R^r$ and
$\bra \hat H\o id\mid \MM_D^\ell\ket\cong\hat H_R^\ell $ commute
inside $\D(H)$. In the next section we will review the
factorization condition of [RS1] guaranteeing
$\M_R(H)\equiv \D(H)\cong\hat H_R^\ell \o \hat H_R^r$.

\sec{The Monodromy Homomorphism}

Following Propositions 2.1 and 2.2 of [M5]
we now provide what may be called the monodromy homomorphism
$\hat H_R\to H$, which as a linear map has already been
discussed in [RS1].

\bsn
{\bf Proposition 6.1} [M5]
{\sl Let $\mon :\hat H\to H$ be the linear map given by
$$
\mon (\varphi) := (\varphi\o id)(R_{op} R),~~\varphi\in\hat H
$$
Then $\mon$ provides an algebra homomorphism
$\mon :\hat H_R\to H$ satisfying for $a\in H$ and $\varphi\in \hat H$}
\beq\label{6.1}
\mon (a\lef\varphi)=a\1\, \mon (\varphi)\,S(a\2)
\eeq
{\bf Proof:} Putting $\A=H$ in Proposition 2.4 we have to check
\beq\label{6.2}
(R^{31} R^{13}) R^{12}(R^{32} R^{23})=R^{12} (\Delta \o id)(R_{op} R)
\eeq
which is straight  forward from the quasitriangularity of $R$. To prove (6.1)
we use the definition (3.2) to compute
\beanon
\mon(a\lef\varphi) &=&
\bra\varphi\o id\mid (S(a\1)\o\one) R_{op} R(a\2\o \one)\ket\\
&=& \bra\varphi\o id\mid(S(a\1)\o\one)R_{op} R(a\2\o a\3 S(a\4))\ket\\
&=&\bra\varphi\o id\mid (S(a\1)a\2\o a\3)R_{op} R(\one\o S(a\4))\ket\\
&=&\bra\varphi\o id\mid(\one\o a\1)R_{op}R(\one\o S(a\2))\ket\\
&=& a\1\, \mon(\varphi)\,S(a\2)
\eeanon
where in the third line we have used that $R_{op}R$ commutes
with $\Delta(H)$. \qed

\bsn
Proposition 6.1 has implicitely been used in [AS] when studying
representations of the monodromy algebra $\hat H_R$ in terms of
representations of $H$.

Next, since $(\Ad a)\,b := a\1bS(a\2),\ a,b\in H$,
defines a Hopf module action of $H$
on itself, Proposition 6.1 immediately implies that $\mon$ extends to a
homomorphism $\Mon$ of the associated crossed products.

\bsn
{\bf Corollary 6.2} [M5]
{\sl The map $\Mon:\hat H\o H\to H\o H$
\beq\label{6.3}
\Mon (\varphi\o a) := \mon (\varphi) \o a
\eeq
provides a homomorphism of algebras $\Mon :\M_R(H)\to
H\cros_{Ad} H.$}

\bsn%
We now note a simple ``bosonization formula''
\footnote{ This terminology is not meant to indicate physical
  interpretations here. It just memorizes the fact, that the
  two non-commuting copies of $H$ generating $H\cros_{Ad} H$ are
  traded for the two commuting copies $U^{-1}(H\o_{alg}\one_H)$
and $U^{-1}(\one_H\o_{alg} H)$.}
showing
that as an algebra the crossed product $H\cros_{Ad} H$ is in fact
isomorphic to $H\o_{alg} H$.

\bsn
{\bf Lemma 6.3}
{\sl Let $U:H\o H\to H\o H$ be given by $U(a\o b):=ab\1\o b\2$.
Then $U$ defines an algebra isomorphism $H\cros_{Ad} H\to H\o_{alg} H$.}

\bsn
{\bf Proof:} $U$ is invertible with $U^{-1} (a\o b)=a S(b\1)\o b\2$. Now the
multiplication in $H\cros_{Ad} H$ is given by
\beq\label{6.4}
(a\o_{Ad} b) (c\o_{Ad} d) = \left(a(\Ad b\1)\,c\o_{Ad} b\2 d\right).
\eeq
Applying $U$ to the r.h.s of (6.4) gives
$$
U(ab\1 c S(b\2)\o_{Ad} b\3 d)=ab\1 cd\1\o b\2 d\2
= [U(a\o_{Ad} b)][U(c\o_{Ad} d)]
$$
where the last product is taken in $H\o_{alg} H$. \qed

\bsn
Note that Lemma 6.3 immediately generalizes to any crossed product by
inner actions.
Putting $\pi_R:= U\circ \Mon$ we now arrive at

\bsn
{\bf Theorem 6.4}
{\sl Let $\MM^\ell\equiv \MM $
and $\MM^r = \RR(\MM^\ell)^{-1} \RR_{op}$ as in (4.3) and
consider the conditions i)-iii) on an algebra homomorphism
 $\pi_R:\M_R(H)\to H\o_{alg} H$
$$
\begin{array}{rrcl}
i)&\qquad \pi_R\circ i_M &=& \Delta_H\\
ii)& (id_H\o \pi_R)(\MM^\ell) &=& R^{21} R^{12}\\
iii)& (id_H\o \pi_R)(\MM^r) &=& R^{13} R^{31}
\end{array}
$$
Then given i) properties ii) and iii) are equivalent with unique solution
$\pi_R=U\circ \Mon$.
In this case $\pi_R$ is an isomorphism if and only if
$\mon$ is bijective implying
$\pi_R^{-1}(H\o_{alg}\one)=\hat H_R^\ell,\
\pi_R^{-1}(\one\o_{alg}H) =\hat H_R^r$, and therefore
$\M_R(H) = \hat H_R^\ell \o_{alg} \hat H_R^r$.}

\bsn
{\bf Proof:} According to
\no{6.3}, (\ref{6.2}) and Proposition 3.4 $\Mon$ is the unique
homomorphism $\M_R(H)\to H\cros_{Ad}H$ satisfying $\Mon \circ i_M=\one_H
\o id_H$ and $(id_H\o \Mon)(\MM^\ell) =R^{21} R^{12} \in H\o (H\o\one_H)
\subset H\o (H\cros_{Ad} H)$.
Hence $\pi_R=U\circ \Mon$ is the unique homomorphism $\M_R(H)\to H\o_{alg} H$
satisfying i) and ii).
Moreover, given i) conditions ii) and iii) are equivalent, since by (2.2)-(2.4)
$$
(id\o\Delta)(R)(R^{21} R^{12})^{-1} (id\o\Delta)(R_{op}) =R^{13} R^{31}
$$
Finally, by (6.3) $\Mon$ and therefore $\pi_R$ are isomorphisms if and only if
$\mon$ is bijective. In this case ii)  implies $\pi_R(\hat H_R^\ell)=
H\o_{alg}\one_H$ and
$\pi_R(\hat H_R^r) = \one_H\o_{alg} \monop(\hat H)$, where
$$
\monop(\varphi):=(\varphi\o id)(R R_{op})\equiv (id\otimes \varphi)(R_{op} R)
$$
provides the monodromy homomorphism $\hat H_R^r\to H$. Since
for $\varphi,\psi\in\hat H$
$$
\bra\varphi\mid\monop(\psi)\ket = \bra\mon(\varphi)\mid\psi\ket
$$
we conclude that
$\monop$ is bijective iff $\mon$ is bijective,
in which case $\pi_R(\hat H_R^r)=\one_H\o_{alg} H$. \qed

\bsn
Theorem 6.4 is the true reason underlying the observation of
[AFFS] on the equivalence of the categories $\Rep(\M_R(H))$ and
$\Rep(H\o_{alg}H)$ provided the map $\mon$ is invertible. In fact, all
representations of $\M_R(H)$ given in [AFFS] are of the form
$\tau\circ\pi_R,\ \tau\in\Rep(H\o_{alg}H)$.
Note, however, that in our approach no assumptions on semi-simplicity
have been made.

\sec{The Hopf Algebra Structure}

In this section we use the isomorphism $\lambda_R :\D(H)\to\M_R(H)$ of
Corollary 5.3
to induce a homomorphism $\Lambda_R :=\pi_R\circ \lambda_R:\D(H)\to
H\o_{alg} H$. We then show that up to cocycle
equivalence $\Lambda_R$ provides a homomorphism of quasitriangular Hopf
algebras, which becomes an isomorphism if and only if $\mon$ is bijectiv.
In this way we recover the result of [RS1], where the invertibility property
of $\mon$  has been called {\em factorizability}.
The results of this section are also reviewed in Chapter 7 of [M8],
however without using our formalism of generating matrices.

First we recall [Dr1]  that $\D(H)$ always  is a quasitriangular
Hopf algebra with coproduct
$\Delta_D$, antipode $S_D$ and  R-matrix $R_D$ given
for $a\in H$ and $\varphi\in\hat H$ by
\bea
\Delta_D(i_D(a)) &=& i_D(a\1)\o i_D (a\2) \label{7.1}\\
\Delta_D (D(\varphi)) &=&
D(\varphi\2) \o D(\varphi\1)\label{7.2}\\
S_D(i_D(a)) &=& i_D(S(a))\label{7.3}\\
S_D (D(\varphi)) &=& D(\hat S^{-1} (\varphi)) \label{7.4}\\
R_D &=& (i_D\o id_{\D(H)})(\DD)\label{7.5}
\eea
This structure may now immediately
be transported to $\M_R(H)$ via $\lambda_R$ to give
$\Delta_M, S_M$ and $R_M$, the formulae for which however turn out to
look less transparent.
\bea
\Delta_M\circ i_M &=& i_M\circ\Delta\label{7.6} \\
(id\o \Delta_M)(\MM)&=& \RR_{op}^{12}\MM^{13} (\RR_{op}^{12})^{-1}
\MM^{12}\label{7.7}\\
S_M\circ i_M &=& i_M\circ S\\
(id\o S_M)(\MM) &=& (S^{-1}\o id)(\RR_{op}^{-1}\MM\RR_{op})\\
R_M &=& (i_M\o id_{\M_R(H)})(\RR_{op}^{-1}\MM)
\eea
where $\RR_{op}$ and $\MM$ have the same meaning as in \no{5.7}.
As it is explained at the end of
this section these results correct some inconsistencies in [AFFS].

We now skip the intermediate object $\M_R(H)$ and study immediately the
composition $\Lambda_R=\pi_R\circ\lambda_R$. To show that up to cocycle
equivalence $\Lambda_R$ provides a Hopf algebra map we first recall the natural
coproduct $\Delta_{H\o H}$ given on $H\o_{alg} H$ by
$$
\Delta_{H\o H} (a\o b) := (a\1\o b\1)\o (a\2\o b\2)
$$
{\bf Lemma 7.1:}
{\sl Let $T:= R_{23}^{-1} \in (H\o H) \o (H\o H)$
\footnote{From now on lower indices refer to the embedding $H\o H\to H^{\o n},
n\ge 3$, and upper indices -- as before -- refer to the embeddings
$\A\o\A\to\A^{\o3}$ or $H\o\A\to\H\o\A\o\A$, respectively, where
$\A=H\o H$. Also, from now on $\o$ always means $\o_{alg}$.}.
Then $T$ is a left $\Delta_{H\o H}$-cocycle and therefore
\beq\label{7.8}
\delta_{H\o H} := \Ad T \circ \Delta_{H\o H}
\eeq
is also coassociative.}

\bsn
{\bf Proof:} Putting $R^{-1} =\sum_i u^i\o v^i$ we have in $(H\o H)^{\o 3}$
\beanon
T^{12} (\Delta_{H\o H} \o id)(T) &=& R^{-1}_{23}\sum [(\one \o u\1^i ) \o
(\one \o u\2^i) \o (v^i \o \one)]\\
 &=& R_{23}^{-1} R_{45}^{-1} R_{25}^{-1} \\
&=& R_{45}^{-1} R_{23}^{-1} R_{25}^{-1} \\
&=& R_{45}^{-1} \sum [(\one \o u^i)\o (v\1^i \o \one) \o (v\2^i \o \one)] \\
&=& T^{23} (id \o \Delta_{H\o H}) (T)
\eeanon
where we have used (2.2) and (2.3).
Hence $T$ is a left $\Delta_{H\o H}$-cocycle.\qed

\bsn
Replacing $(\Delta,R)$ by $(\Delta_{op},R_{op})$, the coproduct
 $\delta_{H\o H}$ has already appeared in Theorem 2.7 of [RS1]
see also § 7.3 of [M8].
Next, we have

\bsn
{\bf Lemma 7.2}
{\sl Let $\R\in(H\o H)\o (H\o H)$ be given by
\beq\label{7.9}
\R:=R_{41}^{-1} R_{42}^{-1} R_{13} R_{23}
\eeq
Then $\R$ is quasi-triangular with respect to $\delta_{H\o H}$.}

\bsn
{\bf Proof:} Since $R$ and $R_{op}^{-1}$ are quasi-triangular w.r.t.
$\Delta,\ \R':=R_{42}^{-1} R_{13}
\in (H\o H)\o (H\o H)$ is quasitriangular w.r.t.
$\Delta_{H\o H}$. The claim now follows from the twist equivalence
\beq \label{7.10}
\R= T_{op} \R' T^{-1}
\eeq
\qed

\bsn
We are now in the position to prove that $\Lambda_R$ provides a
homomorphism of quasitriangular Hopf algebras
$(\D(H),R_D,\Delta_D)\to(H\o H,\delta_{H\o H},\R)$, which in
fact coincides with the homomorphism given in Theorem 2.9 of [RS1],
see also Chapter 7 of [M8].

\bsn
{\bf Theorem 7.3}
{\sl For any quasitriangular $R\in H\o H$ denote
$\Lambda_R:=\pi_R\circ \lambda_R:\D(H)\to H\o H$. Then
$$
\begin{array}{rrcl}
i)&\Lambda_R\circ i_D &=&\Delta\\
ii)& \qquad(id\o \Lambda_R)(\DD) &=& R_{31}^{-1} R_{12}\\
iii) &\delta_{H\o H} \circ \Lambda_R &=&
(\Lambda_R\o \Lambda_R) \circ \Delta_D\\
iv) &(\Lambda_R\o \Lambda_R) (R_D) &=& \R
\end{array}
$$
Moreover, $\Lambda_R$ is uniquely determined by i) and ii) and $\Lambda_R$
is an isomorphism if and only if $\mon$ is bijective.}

\bsn
{\bf Proof:}
Part i) follows from Theorem 6.4i). To prove ii) we use Theorem
6.4ii), (\ref{5.7}) and (\ref{2.2}) to compute
\beq \label{7.11}
(id\o\Lambda_R) (\DD) =
(id\o\Delta)(R_{op}^{-1}) R_{21} R_{12} = R_{31}^{-1}R_{12}
\eeq
The uniqueness of $\Lambda_R$ under conditions i) and ii)
follows from Lemma 5.2.
Next, using i), the claim iii) holds on $i_D(H)\subset \D (H)$ since
$$
\delta_{H\o H}(\Delta(a))
= R_{23}^{-1} (a\1 \o a\3 \o a\2 \o a\4)R_{23}
= \Delta (a\1) \o \Delta (a\2)
$$
It remains to check that iii) holds on all $D(\varphi),\
\varphi\in\hat H$, where $\Delta_D$ is given by $\hat\Delta_{op}$,
according to (\ref{7.2}).
Hence, by (\ref{2.11}) we have to show
\beq\label{7.12}
\FH^{13}\FH^{12} = (id\o\delta_\HH)(\FH) \in H\o(\HH)^{\o 2}
\eeq
where $\FH = (id\o\Lambda_R)(\DD)\in H\o(\HH)$.
Using (\ref{7.11}) we have
$$
\FH^{13}\FH^{12} = R_{51}^{-1}R_{14}R_{31}^{-1}R_{12}
$$
and
\beanon
(id\o\delta_\HH)(\FH)
&=& R_{34}^{-1}\sum(v^ix^j\o y^j\1\o u^i\1\o y^j\2\o u^i\2)R_{34}\\
&=& R_{34}^{-1}R_{51}^{-1}R_{31}^{-1}R_{14}R_{12}R_{34}\\
&=& R_{51}^{-1}(R_{34}^{-1}R_{31}^{-1}R_{14}R_{34})R_{12}\\
&=& R_{51}^{-1}R_{14}R_{31}^{-1}R_{12}\\
&=& \FH^{13}\FH^{12}
\eeanon
where in the second line we have used (\ref{2.2}), (\ref{2.3}) and in
the fourth line the Yang-Baxter equations in the form
$$
R_{14}R_{34}R_{31} = R_{31}R_{34}R_{14}
$$
which follow from (\ref{2.2}) - (\ref{2.5}).
This proves (\ref{7.12}) and therefore iii).
To prove iv) we use (\ref{7.5}), i) and ii) to get
\beanon
(\Lambda_R\o\Lambda_R)(R_D) &=& \Delta(e_\nu)\o (e^\nu\o id\o id)
(R_{31}^{-1}R_{12})\\
&=& (\Delta\o id\o id)(R_{31}^{-1}R_{12})\\
&=& R_{41}^{-1}R_{42}^{-1}R_{13}R_{23}\\
&=& \R
\eeanon
\qed

\bigskip
Let me close by pointing out that [AFFS] seem to provide a
coproduct $\Delta_\K$ on $\K\equiv\M_R(H)$ such that
$$
\Delta'_\HH := (\pi_R\o\pi_R)\circ\Delta_\K\circ\pi_R^{-1}
$$
satisfies
\beq\label{7.13}
\Delta'_\HH(a\o\one) = (a\1\o\one)\o(a\2\o\one)\ .
\eeq
Indeed, assuming $\pi_R$ invertible and defining as in [AFFS]
\bea
\MM_+ &:=& (id_H\o\pi_R^{-1})(R_{op}\o\one_H)\in H\o\M_R(H)\\
\MM_- &:=& (id_H\o\pi_R^{-1})(R^{-1}\o\one_H)\in H\o\M_R(H)
\eea
equ. \no{7.13} would be equivalent to
\beq
(id_H\o\Delta_\K)(\MM_\pm) = \MM_\pm^{12}\MM_\pm^{13}
\eeq
or
\beq
(id_H\o\Delta_\K)(\MM) = \MM_+^{12}\MM^{13}(\MM_-^{12})^{-1}
\eeq
which are the formulas given in Section 4 of [AFFS]. On the
other hand, equ. (4.15) of [AFFS] implies
\beq\label{7.18}
\Delta'_\HH(a\1\o a\2) = (a\1\o a\2)\o(a\3\o a\4)
\eeq
which would be consistent with our $\delta_\HH$. However, equ.
\no{7.13} is manifestly inconsistent with \no{7.18}. In fact,
writing
$$
(\one\o a) = (S(a\1)\o\one)(a\2\o a\3)
$$
it is easy to see that a map $\Delta'_\HH :\HH\to(\HH)\o(\HH)$
obeying \no{7.13} and \no{7.18} cannot consistently be extended
to an algebra homomorphism.

\bsn
After presenting these results
I have been informed [S] that there will be a revised version of [AFFS]
reproducing $\delta_\HH$ -- or equivalently $\Delta_M$
given in \no{7.6},\no{7.7} -- at least up to cocycle
equivalence.

\begin{appendix}

\renewcommand{\theequation}{\mbox{\Alph{section}.\arabic{equation}}}
\sec{Multi-Loop Algebras}

Having identified the monodromy algebra $\hat H_R$ as a braided group
we show in this Appendix that the multi-loop algebras of [AS] arise as braided
tensor products of braided groups in the sense of Majid [M4,M7].
We recall that
these are defined in the natural way so as to obtain all structural maps as
$H$-module morphisms with respect to the coadjoint action $\lef$.
More generally one has\\
\\
{\bf Definition A.1} [M4, M7] Let $(H,\Delta, R=\sum x^i\o y^i)$
be a quasitriangular Hopf
algebra and let $\A$ and $\B$ be left $H$-module algebras with both actions
denoted by $\lef$. The braided tensor product algebra $\A\otimes_R\B$ is
defined to be the vector space $\A\otimes\B$ with multiplication given for
$a,a'\in \A$ and $b,b'\in\B$ by
\beq\label{A.1}
(a\otimes_Rb) (a'\otimes_Rb') := \sum a(y^i \lef a') \otimes_R (x^i \lef b)b'
\eeq

\bsn
One immediately checks that $\A\otimes_R\B$ is a again an $H$-module algebra
with action
\beq\label{A.2}
h\lef (a\otimes_R b) := (h\1 \lef a) \otimes_R (h\2 \lef b)
\eeq
Moreover, the braided tensor product is associative in the sense that
$(\A_1\otimes_R \A_2)\otimes_R \A_3$ and $\A_1 \otimes_R(\A_2 \otimes_R \A_3)$
define the same algebra structure on $\A_1\otimes \A_2\otimes \A_3$.

We now apply this to the multi-loop algebras introduced in Section 6
of [AS] and show that with respect to the coadjoint action they are indeed the
$m$-fold braided tensor product of the associated
one-loop ($\equiv$ monodromy) algebras.
To this end let us denote
\beq\label{A.3}
\L_m := \hat H_R \otimes_R \cdots \otimes_R \hat H_R
\eeq
the $m$-fold braided tensor product and let $\MM_\nu ,\ 1\le
\nu \le m$, denote the $m$ copies of monodromy matrices in $H\otimes \L_m$
$$
\MM_\nu =(id_H\otimes \iota_\nu)(\MM)
$$
where $\iota_\nu :\hat H_R\to \L_m$ is the obvious embedding into the $\nu$-th
tensor factor. We then have\\
\\
{\bf Proposition A.2}
{\sl For $\mu<\nu$ the
following relations hold in $H\otimes
H\otimes \L_m$ }
\beq\label{A.4}
\MM_\nu^{13} R^{12}
\MM_\mu^{23} = R^{12} \MM_\mu^{23} (R^{-1})^{12} \MM_\nu^{13}R^{12}
\eeq
{\bf Proof:} We only proof the case $m=2$, from which the general case follows
straight forwardly. Putting $M_\nu(\varphi):= (\varphi\otimes
id)(\MM_\nu),\ \varphi
\in\hat H$, we conclude from \no{A.1}
$$
M_2(\psi) M_1(\xi) =\sum_i M_1 (y^i\lef \xi)M_2(x^i\lef \psi)
$$
which by \no{3.2} implies
\bea \label{A.5}
\MM_2^{13} \MM_1^{23} &=& \sum_i [S(x^i\1)\otimes S(y\1^i)\otimes \one_\A]
\MM_1^{23} \MM_2^{13} [x\2^i \otimes y^i\2 \otimes \one_\A] \nonumber \\
&=& \sum_{i,j,k,\ell} [S(x^i x^j) \otimes S(y^j y^\ell) \otimes \one_\A]
\MM_1^{23}\MM_2^{13} [x^k x^\ell \otimes y^iy^k \otimes \one_A] \nonumber \\
&=& \sum_\ell [\one_H \otimes S(y^\ell) \otimes \one_\A]
R^{12}\MM_1^{23} (R^{-1})^{12}
\MM_2^{13} R^{12} [x^\ell \otimes \one_H\otimes \one_\A]
\eea
from which the claim follows as in the proof of Proposition 2.4. Here we have
used \no{2.2} and \no{2.3} in the second line and the identities $(S\otimes
id)(R)=R^{-1}$ and $(S\otimes S)(R)=R$ in the last line. \qed

\bsn
A comparison with equ.(6.1) of [AS] shows that up to an ordering convention our
definition of $\L_m$ coincides with their multiloop algebras.

Next we note that the monodromy homomorphism $\mon$ naturally extends to
multi-loop algebras. As a matter of fact, being an $H$-module map intertwining
the coadjoint action on $\hat H_R$ with the adjoint action on $H$ according to
Proposition 6.1, we can use its ordinary $m$-fold tensor product
$
\mon^{\otimes m} = \mon \otimes \cdots \otimes \mon
$
to obtain an $H$-module algebra morphism
\beq\label{A.6}
\mon^{\otimes m} :\L_m\to
H\otimes_R\cdots\otimes_R H
\eeq
where on the r.h.s. the braided tensor product with respect to the left adjoint
action of $H$ on itself is understood.

Next,
similarly as crossed products with respect to inner actions are isomorphic
to algebraic tensor products (see e.g. Lemma 6.3), we now have the
following ``bosonization formula''.\\
\\
{\bf Proposition A.3}
{\sl Consider $H$ as a
left $H$-module algebra under the adjoint
action and
let $\iota_\A :H\to \A$ be a unital algebra map inducing an inner left
$H$-action $\lef$ on $\A$ by
\beq\label{A.7}
h\lef a := \iota_\A(h\1)\, a\,\iota_\A (S(h\2)),\ h\in H,\, a\in \A\ .
\eeq
Then the linear map $V_\A:\A\otimes_R H\to \A\otimes_{alg} H$
\beq\label{A.8}
V_\A (a\otimes_R h):= \sum_{i,j} a\,\iota_\A(y^i y^j) \otimes_{alg}
 x^i h S(x^j)
\eeq
defines an algebra isomorphism
satisfying for all $a\in\A$ and $h,h'\in H$
\beq\label{A.9}
V_\A (h\lef (a\otimes_R h')) =\Delta_\A (h\1) V_\A (a\otimes_R h') \Delta_\A
(S(h\2))
\eeq
where $\Delta_\A:= (\iota_\A\otimes id_H)\circ \Delta$.}

\bsn
{\bf Proof:}
Using $(S\otimes id)(R) =R^{-1}$ and \no{2.14} the inverse of $V_\A$ is
given by
\beq\label{A.10}
V_\A^{-1} (a\otimes_{alg} h)
=\sum_{i,j} a\, \iota_\A(y^j S(y^i)) \otimes_R x^i h x^j
\eeq
Now $V_\A(a\otimes_R\one_H) =(a\otimes_{alg}\one_H)$ and
$$
V_\A(\one_\A\otimes_R h)
=(\iota_\A\otimes id) (R_{op} (\one_H\otimes_{alg} h) R_{op}^{-1}).
$$
proving that the restrictions of $V_\A$ to the subalgebras $\A\otimes_R\one_H$
and $\one_\A\otimes_R H$ are algebra maps. We are left to check that $V_\A$
respects the commutation relations between
$\A\otimes_R\one_H$ and $\one_\A\otimes_R H$.
To this end we compute (implying a summation convention over doubled indices)
\beanon
V_\A((\one_\A\otimes_R h)(a\otimes_R\one_H))
&=& V_\A(y^i \lef a\otimes_R x^i\lef h)\\
&=& y^i\1 a S(y\2^i)y^my^n\, \otimes_{alg}\, x^m x\1^i h S(x\2^i)S(x^n)\\
&=& y^j y^\ell a S(y^iy^k)y^m y^n\,
\otimes_{alg}\, x^m x^i x^j h S(x^n x^kx^\ell)\\
&=& y^jy^\ell a\, \otimes_{alg}\, x^j hS(x^\ell) \\
&=& V_\A (\one_\A\otimes_Rh)\,V_\A (a\otimes_R\one_H)
\eeanon
where we have used \no{2.2} and \no{2.3} in the third line and \no{2.14} in the
fourth line, and where by a
convient abuse of notation we have dropped the symbol
$\iota_\A$.
This proves that $V_\A :\A\otimes_R H\to \A\otimes_{alg} H$ provides
an algebra map. To prove \no{A.9} we use that the general definition \no{A.2}
provides a Hopf module action. Hence, it is enough to check \no{A.9} separately
on the generating factors $\A\otimes_R\one_H\cong \A$
and $\one_\A\otimes_R H\cong
H$. Now, dropping again the symbol $\iota_\A$ we have in $\A\otimes_{alg} H$
\beanon
\Delta_\A(h\1)(a\otimes_{alg}\one_H)\Delta_\A(S(h\2))&=&
h\1 aS(h\4)\otimes_{alg} h\2 S(h\3)\\
&=&(h\lef a)\otimes_{alg}\one_H
\eeanon
proving \no{A.9} on $(\A\otimes_R\one_H)$. On $(\one_\A\otimes_R H)$ we get
\beanon
V_\A(\one_\A \otimes_R h\lef h') &=& y^iy^j \otimes_{alg}
 x^i h\1 h' S(h\2) S(x^j)\\
&=& y^i y^j h\3 S(h\4)\otimes_{alg} x^i h\1 h' S(x^j h\2)\\
&=& h\1 y^i y^j S(h\4) \otimes_{alg} h\2 x^i h' S(x^j) S(h\3)\\
&=& \Delta_\A (h\1) V_\A (\one_\A \otimes_R h') \Delta_\A(S(h\2))
\eeanon
This proves \no{A.9} and therefore Proposition A.3. \qed\\
\\
{\bf Corollary A.4}
{\sl Under the conditions of Proposition A.3 put
$\delta_\A:=V_\A^{-1} \circ \Delta_\A:H\to \A\otimes_R H$. Then $\delta_\A$
implements the $H$-action on $\A\otimes_R H$, i.e.}
$$h\lef (a\otimes_R h') =\delta_\A (h\1) (a\otimes_R h') \delta_\A(S(h\2))$$

\bsn
Note that we may in particular put $\A=H$ and $\iota_\A = id$ to
obtain, using $(S\o S)(R)=R$
\bea\label{A.11}
\delta_\A(a) &=& \sum_{i,j} y^j a\2 y^i \o S^{-1}(x^i)x^j a\1\nonumber\\
&=& \sum_{i,j} a\1 y^jS(y^i)\o x^i a\2 x^j
\eea
which coincides with $\Delta_R$ given in \no{2.1} up to a change of
conventions (i.e. replacing $(H,\Delta,S,R)$ by
$(H,\Delta_{op},S^{-1},R_{op})$.
This shows that $\Delta_R$ provides an algebra map $H\to
H_{cop}\o_{R_{op}}H_{cop}$ as remarked after Lemma 2.2.

We now put
$\iota_{\A\otimes_R H} = \delta_\A$ and proceed inductively to
get isomorphisms
$$V_{\A,m} :\A\otimes_R H\otimes_R\cdots \otimes_R H\to \A\otimes_{alg} H
\otimes_{alg} \cdots \otimes_{alg} H$$ where the tensor factors $H$
appear $m$-times.
Choosing in particular
$\A=\CC$ and $\iota_\A=\e$ and denoting $V_{\CC,m} \equiv
V_m$ we have proven\\
\\
{\bf Theorem A.5}
{\sl The map
$$\mon_{,m} := V_m\circ \mon^{\otimes m} : \L_m\to H^{\otimes m} $$
provides a homomorphism of algebras such that
$$\mon_{,m} (h \lef\underline a) =\Delta^{(m)} (h\1) \mon_{,m}
(\underline a) \Delta^{(m)} (S(h\2))$$
where $\underline a\in \L_m,\ h\in H$
and $\Delta^{(m)} :H\to H^{\otimes m}$ denotes
the $m$-fold coproduct
\footnote{i.e. $\Delta^{(0)} =\e, \Delta^{(1)} =id_H$
and $\Delta^{(m+1)} = (id^{\otimes(m-1)} \otimes \Delta)\circ \Delta^{(m)}$}.
Moreover,
$\mon_{,m}$ is an isomorphism if and only if $\mon$ is bijectiv.}

\bsn
Generalizing Lemma 6.3 in the obvious way we further have
$H^{\o m}\cros_{Ad} H\cong H^{\o (m+1)}$
and therefore also a homomorphism
$$\Mon_{,m}:\L_m\cros H\to H^{\o(m+1)}\,,$$
which is bijective iff $\mon$ is bijective.

Theorem A.5  explains
the representation theory of $\L_m$ given by
[AS] without having to rely on any semi-simplicity assumptions.
To see this explicitly we show\\
\\
{\bf Proposition A.6}
{\sl Let $\MM_\nu
=(id_H\otimes \iota_\nu)(\MM)\in H\otimes \L_m,\
1\le \nu \le m$, be the $\nu$-th monodromy matrix. Then
\beq\label{A.12}
(id_H\otimes \mon_{,m})(\MM_\nu) =\NN_\nu \otimes \one_H^{\otimes(m-\nu)}
\in H\otimes H^{\otimes m}
\eeq
where}
\beq\label{A.13}
\NN_\nu:= (id_H\otimes \Delta^{(\nu-1)} \otimes id_H)(R^{32} R^{31} R^{13}
(R^{-1})^{32})\in H\otimes H^{\otimes \nu}
\eeq
{\bf Proof:} We proceed by induction over $m\ge 1$. The case $m=1$ holds by
definition of $\mon$. Now suppose the
claim holds for $m_0$ and all $1\le\nu \le m_0$.
Putting $\A=H\otimes_R\cdots \otimes_R H\ (m_0$ factors) we have
\beq\label{A.14}
\mon_{,m_0+1} =(V_{m_0}
\otimes id_H)\circ V_\A\circ \mon^{\otimes(m_0+1)}
\eeq
Since the restriction of $V_\A$ for $\A\otimes_R\one_H$ is the identity we get
with respect to the identification $\L_m\cong \L_m\otimes_R\hat\one \subset
\L_{m+1}$
\beq\label{A.15}
\mon_{,m_0+1} | \L_{m_0} = \mon_{,m_0}
\eeq
This proves \no{A.12} for $m=m_0+1$ and $1\le \nu\le m_0$.
For $\nu=m_0+1$ we use
\beq\label{A.16}
(id_H \otimes \mon^{\otimes(m_0+1)} )(\MM_\nu)=R^{31} R^{13}
\in H\otimes_{alg} (\A\otimes_{alg} H)
\eeq
Applying $(id_H\otimes V_\A)$ to
\no{A.16} we get in $H\otimes_{alg}(\A\otimes_{alg} H)$
\beq\label{A.17}
(id_H\otimes V_\A)(R^{31} R^{13})=
(id_H\otimes \iota_\A\otimes id_H)(R^{32} R^{31}
R^{13} (R^{32})^{-1})
\eeq
Finally, according to Corollary A.4
and the inductive definition of $V_m$ we have
$V_m\circ \iota_\A=\Delta^{(m)}$,
which together with \no{A.14}, \no{A.16} and \no{A.17} proves
\no{A.12} for $m=m_0$ and $\nu=m_0+1$.
\qed

\bsn
Comparing \no{A.12} with the representations of $\L_m$ given in Section 6 of
[AS] we realize that they are in fact
all of the form $\tau\circ \mon_{,m}$,
where $\tau\in\Rep (H^{\otimes m})$.

\end{appendix}

\end{document}